\documentclass[11pt,reqno]{amsart}

\usepackage{epsfig}

\usepackage{graphicx}
\usepackage{amscd}
\usepackage{amsmath}
\usepackage{amsfonts}
\usepackage{yfonts}
\usepackage{psfrag}
\usepackage{amssymb}
\usepackage{verbatim}
\usepackage{comment}

\def\bp{{\bf p}}
\newcommand{\Prob}{{\mathbb P}\,}

\def\Pmu{\Prob_{\!\!\mu}}

\newcommand{\be}{\begin{eqnarray}}
\newcommand{\ee}{\end{eqnarray}}

\newcommand{\half}{\frac{1}{2}}

\newcommand{\eps}{{\varepsilon}}

\newcommand{\R}{{\mathbb R}}

\newcommand{\Nat}{{\mathbb N}}
\def\N{\Nat}

\def\Sig{\Sigma}

\newcommand{\gam}{\gamma}

 \newtheorem{theorem}{Theorem}
 \newtheorem{lemma}[theorem]{Lemma}

\input epsf.sty

\begin{document}

\thispagestyle{empty}

\title[Dimension spectrum for a nonconventional average]
{{Dimension spectrum for a nonconventional ergodic average}}

\author{Yuval Peres}
\address{Yuval Peres, One Microsoft Way, Redmond, WA 98052, USA}
\email{peres@microsoft.com}

\author{Boris Solomyak}
\address{Boris Solomyak, Box 354350, Department of Mathematics,
University of Washington, Seattle WA 98195, USA}
\email{solomyak@math.washington.edu}

\date{\today}

\thanks{The research of B. S. was supported in part by the NSF grant DMS-0968879.}


\begin{abstract}
We compute the dimension spectrum of certain nonconventional averages, namely, the Hausdorff dimension of the set of $0,1$
sequences, for which the frequency of the pattern 11 in positions $k, 2k$ equals a given number $\theta\in [0,1]$.
\end{abstract}

\maketitle


\section{Introduction.}

For a dynamical system $(X,T)$ (say, a continuous self-map of a compact metric space), the dimension spectrum of ordinary Birkhoff averages is defined as the function
$$
\theta \mapsto \dim_H\Bigl\{x\in X:\ \lim_{n\to \infty} \frac{S_n f(x)}{n} = \theta \Bigr\}.
$$
where $S_n f(x) = \sum_{k=1}^{n} f(T^k x)$ and $f$ is a function on $X$. It has been widely investigated in Multifractal Analysis, see e.g.\ \cite{barreira_book}.
The most basic example of such analysis goes back to Besicovitch \cite{besic} and Eggleston \cite{egg} who proved that
\be \label{egg}
\dim_H \Bigl\{{(x_k)}_1^\infty\in \{0,1\}^\N:
\ \lim_{n\to \infty} \frac{1}{n} \sum_{k=1}^n x_k = \theta\Bigr\}=H(\theta),\ \ \theta\in [0,1],
\ee
where $H(\theta)= -\theta \log_2\theta - (1-\theta)
\log_2(1-\theta)
$ is the entropy function. Throughout the paper, $\{0,1\}^\N=\Sig_2$ is the symbolic space, with the usual metric
$
\varrho((x_k) ,(y_k)) = 2^{-\min\{n:\ x_n \ne y_n\}}.
$
For dimension purposes, this is equivalent to $[0,1]$ with the standard metric, since for
any set $A\subset \Sigma_2$, its image under the binary representation map has the same dimension as $A$, see \cite[Section 2.4]{falconer}.

Furstenberg \cite{furst} was the first to consider {\em multiple Birkhoff averages}, and their study has become a very active area of research, see e.g.\ Bourgain \cite{bourgain}, Host and Kra \cite{host_kra}, and others.
For a system $(X,T)$ one considers
$$
\frac{1}{n} S_n(f_1,\ldots,f_\ell)(x):=\frac{1}{n} \sum_{k=1}^n f_1(T^k x) f_2(T^{2k}x) \cdots f_\ell(T^{\ell k} x)
$$
for some bounded functions $f_1,\ldots,f_\ell$. Very recently, Yu.\ Kifer \cite{kifer} and A.-H. Fan, L. Liao, J. Ma \cite{FLM} initiated the study of the
dimension spectrum for such averages (in \cite{kifer} more general ``nonconventional averages" are considered as well).
Multifractal analysis of this kind appears to be very complicated, so it is natural to start with the simplest situation, namely, the shift map $T$ on the symbolic space and the functions $f_1,\ldots,f_\ell$ depending only on the first digit $x_1$, for $\ell\ge 2$. Specializing even further, to $\ell=2$ and $f_1(x)\equiv f_2(x)=x_1$ leads to the sets
\be \label{Atheta}
A_\theta:= \Bigl\{{(x_k)}_1^\infty\  \in \Sig_2:
\ \lim_{n\to \infty} \frac{1}{n} \sum_{k=1}^n x_k x_{2k} = \theta\Bigr\}\,,\ \theta\in [0,1].
\ee
The question about the dimension of $A_\theta$ was raised in \cite{FLM}. Note that this directly generalizes the Besicovitch-Eggleston set-up from $\ell=1$ to $\ell=2$.

Motivated by this problem, A.-H. Fan, L. Liao, J. Ma, and J. Schmeling [private communication in August 2010]
computed the Minkowski (box-counting) dimension of another set
$$
X_G:= \Bigl\{ {(x_k)}_1^\infty\  \in \Sig_2:\ x_k x_{2k}=0 \ \mbox{for all}\ k\Bigr\}
$$
and asked what is its Hausdorff dimension.
It is obvious that $X_G\subset A_0$, and in fact, it is easy to see that $\dim_H(X_G)=\dim_H(A_0)$.

In joint work with R. Kenyon, we computed the Hausdorff dimension of $X_G$ and a large class of
similarly defined sets, putting it into the context of {\em subshifts invariant under the semi-group of multiplicative integers}
\cite{KPS1,KPS2}. Here we adapt the techniques of \cite{KPS1,KPS2} to compute the full
dimension spectrum $\dim_H(A_\theta)$.

\begin{theorem} \label{th-freq}
Let $A_\theta$ be given by (\ref{Atheta}). For $\theta\in (0,1)$ we have
\be \label{eq-dim1}
\dim_H(A_\theta) = f(\theta):=-\log_2(1-p)-\frac{\theta}{2} \log_2\Bigl[\frac{(1-q)(1-p)}{qp}\Bigr],
\ee
where
\be \label{eq-pq}
p^2q = (1-p)^3, \ 0 <p < 1, \ 0 < q < 1,
\ee
\be \label{eq-theta}
\theta = \frac{2p(1-q)}{1+p+q}\,.
\ee
We have $\dim_H(A_0)=\lim_{\theta\to 0} f(\theta)=-\log_2(1-p)$, with $p^2=(1-p)^3$, and $\dim_H(A_1)=\lim_{\theta\to 1} f(\theta)=0$.
\end{theorem}

The meaning of $p$ and $q$ will be explained in the next section. Of course,
it is easy to eliminate $q$ from (\ref{eq-dim1}) and (\ref{eq-theta}). For a given $\theta$, we get an algebraic equation of degree 4 for
$p$. Solving the equation numerically yields the graph in Figure 1.

\begin{figure}[ht]
\includegraphics[height=2.4in]{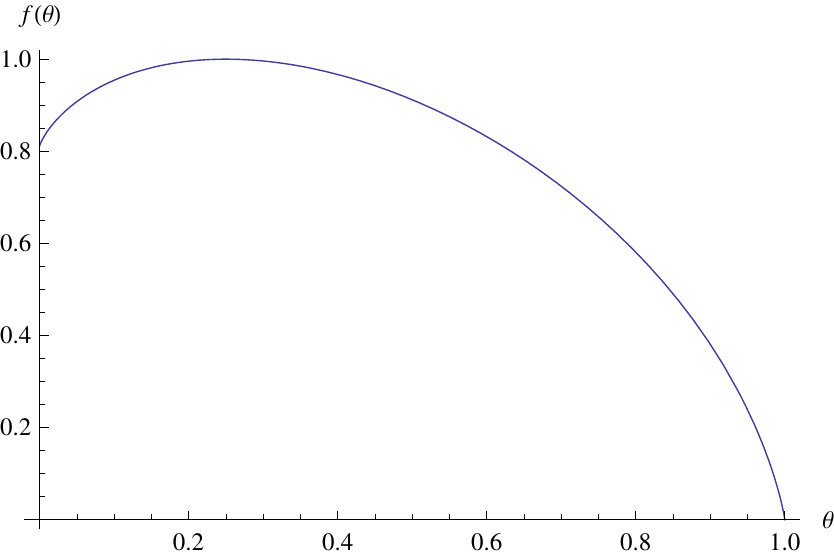}
\caption{Dimension of $A_\theta$}
\end{figure}


\medskip

\noindent {\bf Remarks.}

1. As already mentioned, the formula for $\dim_H(A_0)$ easily follows from \cite{KPS1}. Note that in
\cite{KPS1,KPS2} notation was slightly different, so that $p$ in those papers is $1-p$ here.

2. It is immediate that $A_1$ is contained in the set of $0$-$1$ 
sequences which have frequency of 1's equal to 1. Thus,
$\dim_H(A_1)=0$ by (\ref{egg}), and we assume $\theta<1$ for the rest of the paper.

\begin{sloppypar}
3. By the Strong Law of Large Numbers for weakly correlated random variables (see \cite{lyons}), for a.e.\ sequence $(x_k)$
with respect to the Bernoulli $(\half,\half)^\N$
measure, $\lim_{n\to \infty} \frac{1}{n} \sum_{j=1}^n x_j x_{2j}=1/4$. This agrees with our result: for $\theta=1/4$ we get $\dim_H(A_\theta)=1$ and $p=q=1/2$.
\end{sloppypar}


4. In \cite{FLM} it is proved that
$\dim_H(B_\theta) = 1-\frac{1}{\ell} + \frac{1}{\ell} H\left(\frac{1+\theta}{2}\right)$ for $\theta\in [-1,1]$ and $\ell\ge 2$,
where
$B_\theta:= \bigl\{{(x_k)}_1^\infty\ \in \{-1,1\}^\N:
\ \lim_{n\to \infty} \frac{1}{n} \sum_{k=1}^n x_k x_{2k}\cdots x_{\ell k} = \theta\bigr\}$,
using the techniques of Riesz products. It is further pointed out in \cite{FLM}
that the problem becomes drastically different if one takes the digits 0,1 (which reduces to $A_\theta$ for $\ell=2$)
instead of $-1,1$.

5. Yu.\ Kifer \cite{kifer} considered a slightly different question: he studied the Hausdorff dimension of sets defined by
the frequencies of all $\ell$-tuples of digits $i_1,\ldots,i_\ell$ in positions $k,2k,\ldots, \ell k$. However, he was able to compute the dimensions only under
the assumption that such frequencies are of the form $p_{i_1,\ldots,i_\ell}=p_{i_1}\cdots  p_{i_\ell}$.

6. As pointed out in \cite{KPS2}, there are some parallels between the multiplicative shifts of finite type and
self-affine carpets \cite{Bedf,McM}; we should also add here self-affine sponges \cite{KP}. The present paper may similarly be compared to the work on
multifractal self-affine carpets and sponges, see e.g.\ \cite{king,olsen,BaMe,JR}; however, we do not see any way to directly
transfer the results.


\section{Preliminaries and the scheme of the proof.}
The dimension of $A_\theta$ is computed with the help of the following lemma which goes back to Billingsley \cite{Billing}.
We write $[u]$ for the cylinder set of sequences starting with
a finite word $u$ and $x_1^n: = x_1\ldots x_n$.

\begin{lemma}[see Prop.4.9 in \cite{falconer}] \label{prop-mass}
Let $E$ be a Borel set in $\Sig_2$ and let $\nu$ be a finite Borel measure on $\Sig_2$.

{\bf (i)} If $\nu(E)>0$ and
$\liminf_{n\to \infty} \frac{-\log_2 \nu[x_1^n]}{n} \ge s\ \ \mbox{for $\nu$-a.e.}\ x\in E,$
then $\dim_H(E) \ge s$.

{\bf (ii)} If
$\liminf_{n\to \infty} \frac{- \log_2 \nu[x_1^n]}{n} \le s\ \ \mbox{for all}\ x\in E,$
then $\dim_H(E) \le s$.
\end{lemma}

Following \cite{KPS1,KPS2},
for a probability measure $\mu$ on $\Sig_2$, we define another measure $\Pmu$ on $\Sig_2$ by
\be \label{eq-meas1}
\Pmu[u]:= \prod_{i\le n,\, i\ \mbox{\tiny odd}} \mu[u|_{J(i)}], \ \ \mbox{where}\ J(i) = \{2^r i\}_{r=0}^\infty
\ee
and $u|_{J(i)}$ is the subsequence of $u$ (viewed as a finite sequence) along the geometric progression $J(i)$.
The new measure $\Pmu$ is invariant under the action of the {\em multiplicative semigroup of odd positive numbers}:
$$
(x_k)_{k=1}^\infty\mapsto (x_{ik})_{k=1}^\infty\ \ \mbox{for odd}\ i.
$$

We consider Markov measures $\mu_{\bp,P}$ on $\Sig_2$, with the initial probability distribution $\bp=(1-p,p)$
(so that $p$ is the probability of initial 1), and the stochastic transition
matrix $P = \left( \begin{array}{cc} 1-p & p \\ q & 1-q \end{array} \right)$.
Note that our Markov measures are not stationary; instead, their initial distribution coincides with the first row of the transition matrix.

Next we indicate the scheme of the proof of Theorem~\ref{Atheta}. Recall that
$\theta\in [0,1)$.
In view of Lemma~\ref{prop-mass}(i),
the lower bound for $\dim_H(A_\theta)$ will be established once we  prove the following.

\begin{lemma} \label{prop-Markov}
Fix $p\in (0,1), q\in [0,1)$, and let $\Pmu$, with $\mu = \mu_{\bp,P}$, be defined by (\ref{eq-meas1}).

{\bf (i)} If $p,q$ satisfy (\ref{eq-theta}), then $\Pmu(\Sig_2\setminus A_\theta)=0$. For $\theta=0$ we take $q=1$.

{\bf (ii)} For any $p,q$ we have
$$
\lim_{n\to \infty} \frac{-\log_2\Pmu[x_1^n]}{n} = s(p,q):= \frac{(1+q)H(p) + pH(q)}{1+p+q} \ \ \mbox{for $\Pmu$-a.e.}\ x.
$$

{\bf (iii)} The maximum of $s(p,q)$, subject to (\ref{eq-theta}), is achieved when $p^2 q = (1-p)^3$, and it equals 
\be \label{eq-smax}
f(\theta)=-\log_2(1-p)-\frac{\theta}{2} \log_2\Bigl[\frac{(1-q)(1-p)}{qp}\Bigr]
\ee
\end{lemma}

The upper bound in Theorem~\ref{th-freq} will follow from Lemma~\ref{prop-mass}(ii), once we prove the following

\begin{lemma} \label{lem-upper} Let $\mu=\mu_{\bp,P}$ be the Markov measure with initial probability vector
$\bp = (1-p,p)$ and transition matrix $P = \left( \begin{array}{cc} 1-p & p \\ q & 1-q \end{array} \right)$, where
$p^2q = (1-p)^3$ and (\ref{eq-theta}) holds, and let $\Pmu$ be the corresponding multiplicative invariant measure.
Then
\be \label{eq-upper1}
\liminf_{n\to \infty} \frac{-\log_2\Pmu[x_1^n]}{n}  \le  f(\theta)\ \ \mbox{\rm for all}\ x\in A_\theta.
\ee
\end{lemma}


\section{Proofs}

The following elementary lemma will be useful. We provide the proof for completeness.

\begin{lemma}\label{lem-elem}
Suppose that $\{z_n\}$ is a bounded real sequence and there exists $c>0$ such that
\be \label{estik}
|z_n-z_{n+m}|\le \frac{cm}{n}\ \ \mbox{for all}\ m,n\in \N.
\ee
If $z_{2^k n}\to \gamma$ as $k\to \infty$ for all $n\in \N$, then $z_n\to \gamma$.

\end{lemma}

\begin{proof}
For $\eps>0$ let $i_0\in \N$ be such that $2^{-i_0}<\eps$. By the assumption, we can find $k_0\in \N$ such that
for all $\ell\le 2^{i_0}$ and all $k\ge k_0-i_0$ we have $|z_{\ell\cdot 2^k}-\gam|<\eps$. For $n>2^{k_0+i_0}$, with $2^{k-1}\le n < 2^k$,
let $\ell= \lceil n\cdot 2^{i_0-k}\rceil\le 2^{i_0}$. Then $0 \le \ell \cdot 2^{k-i_0}-n< 2^{k-i_0}$, hence
$$
|z_n - \gam| \le |z_n-z_{\ell\cdot 2^{k-i_0}}|+|z_{\ell\cdot 2^{k-i_0}}-\gam|\le
\frac{c2^{k-i_0}}{n}+\eps < \eps(1+2c),
$$
completing the proof.
\end{proof}

For positive integers $m<n$, $i,j\in \{0,1\}$, and $x\in \Sig_2$ let
$$
N_{i}(x_m^n) = \#\left\{k\in [m,n]\cap\N:\ x_k=i\right\},
$$
$$
N_{ij}(x_m^n) = \#\left\{k\in [m,n]\cap\N:\ x_k=i,\ x_{2k}=j\right\}.
$$
Further, for $x\in \Sig_2$ and even $n$ denote
\be \label{def-alpha}
\alpha_n(x) = \frac{N_1(x_{n/2+1}^n)}{n/2}\,.
\ee
Observe that $\alpha_n(x)\in [0,1]$, and it is easy to see that 
\be \label{estik1}
|\alpha_n-\alpha_{n+m}|\le \frac{2m}{n}\ \ \mbox{for} \ n,m \ \mbox{even.}
\ee

\begin{lemma}\label{lem-prep}
Let $\mu$ be  a Markov measure on $\Sig_2$, with initial probability
(row) vector $\bp=(1-p,p)$ and transition matrix $P=\left( \begin{array}{cc} 1-p & p \\ q & 1-q \end{array} \right)$.
Let $\Pmu$ be the measure on $\Sig_2$ defined by (\ref{eq-meas1}). Then for $\Pmu$-a.e.\ $x\in \Sig_2$ we have
$$
\lim_{n\to\infty} \alpha_n(x)=\frac{2p}{1+p+q}=:\xi,
$$
where $\alpha_n(x)$ is defined in (\ref{def-alpha}).
\end{lemma}

\begin{proof}
We assume $n$ to be even.
Consider $n\alpha_{2n}(x)$, the number of 1's in $x_i$'s for $i\in (n,2n]$. 
There are $n/2$ odd numbers $i$ in $(n,2n]$. By the definition of the measure $\Pmu$, 
the corresponding $x_i$ are chosen to be 1 independently
with probability $p$. Thus, by the Law of Large Numbers, the number of 1's among them is $np/2 +o(n)$. Among the even numbers $2i$
in $(n/2,n]$, there is $n\alpha_n(x)/2$ numbers for which $x_i=1$, 
and there will be approximately $n \alpha_n(x)(1-q)/2$ \ 1's there, again using the 
Law of Large Numbers, since $1-q$ is the probability of the transition $1\to 1$ according to the Markov measure $\mu$. Similarly, there is
$n(1-\alpha_n(x))/2$ numbers $2i\in (n,2n]$ for which $x_i=0$, and approximately $n(1-\alpha_n(x))p/2$ 
\ 1's among those. Hence we have
\begin{eqnarray*}
\alpha_{2n}(x) & = & \frac{p}{2} + \frac{\alpha_n(x)(1-q)}{2} + \frac{\alpha_n(x)(1-\alpha_n(x))p}{2} + \eps_n(x)\\
                   & = & p + \frac{\alpha_n(x) (1-q-p)}{2} + \eps_n(x),
\end{eqnarray*}
where $\eps_n(x)\to 0$ for $\Pmu$-a.e.\ $x$. 
Since $|(1-q-p)/2|< 1$, it follows that  
$\alpha_{2^k n}(x)\to \xi$ for $\Pmu$-a.e.\ $x$, where 
$\xi = p+ \xi(1-q-p)/2$, i.e. $\xi = \frac{2p}{1+p+q}$.  Then for $\Pmu$-a.e.\ $x$ we have
$\alpha_{2^kn}(x)\to \xi$ for all $n\in \N$.
It remains to recall (\ref{estik1}) and apply Lemma~\ref{lem-elem}.
\end{proof}

\medskip

\begin{proof}[Proof of Lemma~\ref{prop-Markov}{\rm (i)}]
We claim that 
$$
\lim_{n\to \infty} \frac{N_{11}(x_1^n)}{n}=\xi (1-q)=\frac{2p(1-q)}{1+p+q}=\theta
$$
 for a $\Pmu$-typical point $x\in \Sig_2$.
 Indeed, the initial part of the sequence does not impact the limit, and for $n$ big enough we see 1's with frequency $\xi$ and thus
11's (in positions $i, 2i$) with frequency $\xi(1-q)$, by the definition of $\Pmu$. This means that $\Pmu$-a.e.\ $x$ is in 
$A_\theta$, as desired.
\end{proof}

Fix an even integer $n$. Denote
$$
N_{1,\mbox{\tiny odd}}=N_{1,\mbox{\tiny odd}}(x_1^n):= \#\{k\le n:\ k\ \mbox{odd},\ x_k=1\}.
$$
By the definition of $\mu=\mu_{\bp,P}$ and $\Pmu$ we have, for any $x$ and even $n$:
\be \label{def-Pmu}
\Pmu[x_1^n] = p^{N_{1,\mbox{\tiny odd}}}(1-p)^{n/2 - N_{1,\mbox{\tiny odd}}}
             (1-p)^{N_{00}} p^{N_{01}}
              q^{N_{10}} (1-q)^{N_{11}}
\ee
where $N_{ij}=N_{ij}(x_1^{n/2})$. 

\medskip

For the rest of the proof we write $\log$ for $\log_2$ to simplify notation.

\begin{proof}[Proof of Lemma~\ref{prop-Markov}{\rm (ii)}] 
For $\Pmu$-a.e.\ $x\in \Sig_2$, we see 1's with frequency $p$ in odd places, and similarly to the proof of 
Lemma~\ref{prop-Markov}(i), the frequency of 00's, 01's, 10's, and 11's is $(1-\xi)(1-p)$, $(1-\xi)p$, $\xi q$, and $\xi(1-q)$
respectively. Therefore, for $\Pmu$-a.e.\ $x\in \Sig_2$ the formula (\ref{def-Pmu}) yields
\begin{eqnarray*}
 \lim_{n\to \infty} \frac{-\log \Pmu[x_1^n]}{n}
& = & -(1/2)\bigl[(2-\xi)p\log p + (2-\xi)(1-p)\log(1-p)\\ 
&   & +  \xi q\log q + \xi (1-q)\log(1-q)\bigr]\\[1.2ex]
& = & \frac{(1+q)H(p)+pH(q)}{1+p+q}\,,
\end{eqnarray*}
as desired.
\end{proof}

In view of Lemma~\ref{prop-Markov}(i),(ii) and Lemma~\ref{prop-mass}(i), we have that
\be \label{eq-cons}
\dim_H(A_\theta) \ge s(p,q) = \frac{(1+q)H(p)+pH(q)}{1+p+q},\ \ \mbox{where}\ \
\theta = \frac{2p(1-q)}{1+p+q}\,.
\ee
Thus, we should find the constrained maximum of $s(p,q)$ on $[0,1]^2$. This is a straightforward
exercise, but we include it for the record and in order to explain where the formula (\ref{eq-pq}) comes from.
It is actually not needed for the proof of the main result, since we could just
produce the answer for the optimization problem and refer to Lemma~\ref{lem-upper}. We also include the verification of (\ref{eq-smax}).

\begin{proof}[Proof of Lemma~\ref{prop-Markov}{\rm (iii)}]
We use the method of Lagrange multipliers. Differentiating $s(p,q)$ yields
\begin{eqnarray*}
(1+p+q)^2 \frac{\partial s(p,q)}{\partial p}
& = & (1+q)[(1+p+q)\log(\mbox{$\frac{1-p}{p}$})-H(p)+H(q)]\\
& = & (1+q)[(2+q)\log(1-p)-(1+q)\log p +H(q)],
\end{eqnarray*}
\begin{eqnarray*}
(1+p+q)^2 \frac{\partial s(p,q)}{\partial q}
& = & p[(1+p+q)\log(\mbox{$\frac{1-q}{q}$})+H(p)-H(q)]\\
& = & p[(2+p)\log(1-q)-(1+p)\log q +H(p)].
\end{eqnarray*}
Differentiating the constraint $g(p,q) = \theta(1+p+q)-2p(1-q)=0$ yields
$$
\nabla g(p,q) = (\theta-2(1-q), \theta+2p) =
\Bigl(\frac{-2(1-q)(1+q)}{1+p+q}, \frac{2p(2+p)}{1+p+q}\Bigr).
$$
At the point of constrained maximum we have $\nabla s(p,q) = \lambda \nabla g(p,q)$,
which reduces to
\begin{eqnarray*}
& & (2+p)[(2+q)\log(1-p)-(1+q)\log p +H(q)]\\
& = & -(1-q)[(2+p)\log(1-q)-(1+p)\log q +H(p)].
\end{eqnarray*}
The latter becomes, after collecting the terms:
$$
3(1+p+q)\log(1-p)=2(1+p+q)\log p + (1+p+q)\log q,
$$
so $p^2q=(1-p)^3$, as claimed.

It remains to verify the formula (\ref{eq-smax}). We have
\begin{eqnarray*}
f(\theta) & = & -\log(1-p)-\frac{\theta}{2} \log\Bigl[\frac{(1-q)(1-p)}{qp}\Bigr]\\
          & = & -\log(1-p)- \frac{p(1-q)}{1+p+q} \log\Bigl[\frac{(1-q)(1-p)}{qp}\Bigr].
\end{eqnarray*}
Comparing the latter with
$$
s(p,q) = \frac{(p\log \frac{1-p}{p} - \log(1-p))(1+q) + p(q\log \frac{1-q}{q} - \log(1-q))}{1+p+q}
$$
results in
$
(1+p+q)(f(\theta)-s(p,q))= -p[2\log \frac{1-p}{p} + \log \frac{1-q}{q} + \log(1-p) - \log(1-q)]=0,
$
whenever $p^2q=(1-p)^3$, as desired.
\end{proof}

\begin{proof}[Proof of Lemma~\ref{lem-upper}]
Fix any $x$ and an even $n\in \N$. We continue to use the notation $N_{ij} = N_{ij}(x_1^{n/2})$.
Since $N_{00} + N_{01} + N_{10} + N_{11} = n/2$, we have from (\ref{def-Pmu}):
\begin{eqnarray*}
\Pmu[x_1^n] & = & p^{N_{1,\mbox{\tiny odd}}+N_{01}}(1-p)^{n - N_{1,\mbox{\tiny odd}}-N_{01}-N_{10}-N_{11}}
                q^{N_{10}} (1-q)^{N_{11}} \\
            & = & (1-p)^n {\Bigl(\frac{p}{1-p}\Bigr)}^{N_{1,\mbox{\tiny odd}}+N_{01}}
                          {\Bigl(\frac{q}{1-q}\Bigr)}^{N_{10}}
                          {\Bigl(\frac{1-q}{1-p}\Bigr)}^{N_{10}+N_{11}}.
\end{eqnarray*}
Observe that
$$
N_{10}+N_{11}=N_1(x_1^{n/2})\ \ \ \ \mbox{and}\ \ \ \  N_{1,\mbox{\tiny odd}}=N_1(x_1^n)-N_{01} - N_{11}.
$$
The equation $p^2q=(1-p)^3$ can be rewritten as $\frac{1-q}{1-p}=\frac{1-q}{q}{\bigl(\frac{1-p}{p}\bigr)}^2$. Combining this
with the last several equalities yields
\begin{eqnarray*}
\Pmu[x_1^n] & = & (1-p)^n {\Bigl(\frac{p}{1-p}\Bigr)}^{N_1(x_1^n)-N_{11} - 2N_1(x_1^{n/2})}{\Bigl(\frac{1-q}{q}\Bigr)}^{N_{11}}  \\
            & = & (1-p)^n {\Bigl(\frac{(1-q)(1-p)}{qp}\Bigr)}^{N_{11}} {\Bigl(\frac{p}{1-p}\Bigr)}^{N_1(x_1^n)-2N_1(x_1^{n/2})}.
\end{eqnarray*}
Thus,
\begin{eqnarray}
\frac{-\log\Pmu[x_1^n]}{n} & = & -\log(1-p) - \frac{N_{11}(x_1^{n/2})}{n} \log\Bigl(\frac{(1-q)(1-p)}{qp}\Bigr) \nonumber \\ 
                           &   & + \Bigl(\frac{N_1(x_1^n)}{n} - \frac{N_1(x_1^{n/2})}{n/2}\Bigr) \log\Bigl(\frac{p}{1-p}\Bigr)\nonumber \\
                            & = & f(\theta) + \Bigl(\theta - \frac{N_{11}(x_1^{n/2})}{n}\Bigr)\log\Bigl(\frac{(1-q)(1-p)}{qp}\Bigr) \nonumber \\
                            &   & + \Bigl(\frac{N_1(x_1^n)}{n} - \frac{N_1(x_1^{n/2})}{n/2}\Bigr) \log\Bigl(\frac{p}{1-p}\Bigr), \label{equ2}
\end{eqnarray}
where $f(\theta)$ is from (\ref{eq-smax}).
Observe that $\lim_{n\to \infty} \frac{N_{11}(x_1^{n/2})}{n} = \theta/2$ for all $x\in A_\theta$.
Now replace $n$ by $2^\ell$ for $\ell=1,\ldots,L$, and take the average over $\ell$. 
The expression in the last line of (\ref{equ2}) ``telescopes,'' so we obtain
\begin{eqnarray*}
\frac{1}{L} \sum_{\ell=1}^L \Bigl( \frac{-\log\Pmu[x_1^{2^\ell}]}{2^\ell} - f(\theta) \Bigr)
 & \!= \!& \frac{1}{L} \sum_{\ell=1}^L \Bigl( \frac{\theta}{2}-\frac{N_{11}(x_1^{2^\ell})}{2^\ell} \Bigr) 
\log\Bigl(\frac{(1-q)(1-p)}{qp}\Bigr) \\
 & \! +  \! & \frac{1}{L} \log\Bigl(\frac{p}{1-p}\Bigr) \Bigl(\frac{N_1(x_1^{2^L})}{2^L} - \frac{N_1(x_1^2)}{2}\bigr),
\end{eqnarray*}
which tends to zero, as $L\to \infty$. It follows that
$$
\liminf_{\ell\to \infty} \frac{-\log\Pmu[x_1^{2^\ell}]}{2^\ell} \le f(\theta)\ \ \mbox{for all}\ x\in A_\theta,
$$
and the proof of (\ref{eq-upper1}) is complete.
\end{proof}

\section{Concluding remarks}

1. It is not hard to verify that, under the conditions (\ref{eq-pq}) and (\ref{eq-theta}) we have $\frac{(1-q)(1-p)}{qp}<1$ if and only if
$\theta<1/4$. Therefore, by the argument in the last section, it immediately follows that
$$
\dim_H(A^+_\theta)=f(\theta)\ \ \mbox{for}\ \theta\in (0,1/4),\ \ \ \ 
\dim_H(A^-_\theta)=f(\theta)\ \ \mbox{for}\ \theta\in (1/4,1),
$$
where
$$
A^+_\theta:= \Bigl\{{(x_k)}_1^\infty\  \in \Sig_2:
\ \limsup_{n\to \infty} \frac{1}{n} \sum_{k=1}^n x_k x_{2k} \le \theta\Bigr\},
$$
$$
A^-_\theta:= \Bigl\{{(x_k)}_1^\infty\  \in \Sig_2:
\ \liminf_{n\to \infty} \frac{1}{n} \sum_{k=1}^n x_k x_{2k} \ge \theta\Bigr\}.
$$

2. We extended the result of Theorem~\ref{th-freq} to the case of arbitrary functions $f_1,f_2$ on the shift $\Sig_2$ depending on the first digit $x_1$. 
The method is the same, but the calculations are more involved, so we only state the result.

\begin{theorem} \label{th-freq2}
For $\beta,\gamma\in \R$, let
$$
A_\theta(\beta,\gamma):=
\Bigl\{{(x_k)}_1^\infty\  \in \Sig_2:
\ \lim_{n\to \infty} \frac{1}{n} \sum_{k=1}^n (x_k+\beta) (x_{2k}+\gamma) = \theta\Bigr\}
$$
We have
$$
\dim_H(A_\theta(\beta,\gamma))=
-\half\log_2[p_0(1-p)]-\frac{\theta}{2} \log_2\Bigl[\frac{(1-q)(1-p)}{qp}\Bigr],
$$
where
$$
\frac{1-p}{1-q} = {\Bigl(\frac{q}{1-q}\Bigr)}^{1+2\beta+\gamma} {\Bigl(\frac{p}{1-p}\Bigr)}^{2+2\beta+\gamma}\,,\ \ \ \ 
\frac{1-p_0}{p_0} = {\Bigl(\frac{q}{1-q}\Bigr)}^{\beta} {\Bigl(\frac{p}{1-p}\Bigr)}^{1+\beta}\,,
$$
and
$$
\theta = \beta\gamma + \frac{(1+\beta+\gamma-q)(1+p-p_0) + \beta(p_0(p+q)-q)}{1+p+q}\,.
$$
The appropriate measure is $\Pmu$, with $\mu$ Markov, having the initial distribution $(p_0,1-p_0)$ and the transition matrix 
$\left( \begin{array}{cc} 1-p & p \\ q & 1-q \end{array} \right)$.
\end{theorem}

3. After this work was essentially completed, we were informed that A.-H. Fan, J. Schmeling, and M. Wu have computed the
dimension of $A_\theta$ (in a different, but equivalent form) and other sets of this type, independently, but also building on
\cite{KPS2}.

\medskip

\noindent {\bf Acknowledgment.} We are grateful to the referee for helpful comments and suggestions.



\end{document}